
\frenchspacing

\parindent15pt

\abovedisplayskip4pt plus2pt
\belowdisplayskip4pt plus2pt 
\abovedisplayshortskip2pt plus2pt 
\belowdisplayshortskip2pt plus2pt  

\font\twbf=cmbx10 at12pt
 at12pt
 at12pt

\font\sc=cmcsc10

\font\ninerm=cmr9 
\font\nineit=cmti9 
\font\ninesy=cmsy9 
\font\ninei=cmmi9 
\font\ninebf=cmbx9 

\font\sevenrm=cmr7  
 
\font\seveni=cmmi7  
\font\sevensy=cmsy7 

\font\fivenrm=cmr5  
\font\fiveni=cmmi5  
\font\fivensy=cmsy5 

\def\nine{%
\textfont0=\ninerm \scriptfont0=\sevenrm \scriptscriptfont0=\fivenrm
\textfont1=\ninei \scriptfont1=\seveni \scriptscriptfont1=\fiveni
\textfont2=\ninesy \scriptfont2=\sevensy \scriptscriptfont2=\fivensy
\textfont3=\tenex \scriptfont3=\tenex \scriptscriptfont3=\tenex
\def\rm{\fam0\ninerm}%
\textfont\itfam=\nineit    
\def\it{\fam\itfam\nineit}%
\textfont\bffam=\ninebf 
\def\bf{\fam\bffam\ninebf}%
\normalbaselineskip=11pt
\setbox\strutbox=\hbox{\vrule height8pt depth3pt width0pt}%
\normalbaselines\rm}

\hsize30cc
\vsize44cc
\nopagenumbers

\def\luz#1{\luzno#1?}
\def\luzno#1{\ifx#1?\let\next=\relax\yyy
\else \let\next=\luzno#1\xxx\fi\next}
\def\sp#1{\def\xxx{\kern1.7pt}\def\yyy{\kern-1.7pt}\luz{#1}}
\def\spa#1{\def\xxx{\kern1pt}\def\yyy{\kern-1pt}\luz{#1}}

\newcount\beg
\newbox\aabox
\newbox\atbox
\newbox\fpbox
\def\abbrevauthors#1{\setbox\aabox=\hbox{\sevenrm\uppercase{#1}}}
\def\abbrevtitle#1{\setbox\atbox=\hbox{\sevenrm\uppercase{#1}}}
\long\def\pag{\beg=\pageno
\def\leftheadline{\noindent\rlap{\nine\folio}\hfil\copy\aabox\hfil}
\def\rightheadline{\noindent\hfill\copy\atbox\hfill\llap{\nine\folio}}
\def\phead{\setbox\fpbox=\hbox{\sevenrm 
************************************************}%
\noindent\vbox{\sevenrm\baselineskip9pt\hsize\wd\fpbox%
\centerline{***********************************************}

\centerline{BANACH CENTER PUBLICATIONS, VOLUME **}

\centerline{INSTITUTE OF MATHEMATICS}

\centerline{POLISH ACADEMY OF SCIENCES}

\centerline{WARSZAWA 19**}}\hfill}
\footline{\ifnum\beg=\pageno \hfill\nine[\folio]\hfill\fi}
\headline{\ifnum\beg=\pageno\phead
\else
\ifodd\pageno\rightheadline \else \leftheadline \fi 
\fi}}

\newbox\tbox
\newbox\aubox
\newbox\adbox
\newbox\mathbox

\def\title#1{\setbox\tbox=\hbox{\let\\=\cr 
\baselineskip14pt\vbox{\twbf\tabskip 0pt plus15cc
\halign to\hsize{\hfil\ignorespaces \uppercase{##}\hfil\cr#1\cr}}}}

\newbox\abbox
\setbox\abbox=\vbox{\vglue18pt}

\def\author#1{\setbox\aubox=\hbox{\let\\=\cr 
\nine\baselineskip12pt\vbox{\tabskip 0pt plus15cc
\halign to\hsize{\hfil\ignorespaces \uppercase{\spa{##}}\hfil\cr#1\cr}}}%
\global\setbox\abbox=\vbox{\unvbox\abbox\box\aubox\vskip8pt}}

\def\address#1{\setbox\adbox=\hbox{\let\\=\cr 
\nine\baselineskip12pt\vbox{\it\tabskip 0pt plus15cc
\halign to\hsize{\hfil\ignorespaces {##}\hfil\cr#1\cr}}}%
\global\setbox\abbox=\vbox{\unvbox\abbox\box\adbox\vskip16pt}}

\def\mathclass#1{\setbox\mathbox=\hbox{\footnote{}{1991 {\it Mathematics Subject 
Classification}\/: #1}}}

\long\def\maketitlebcp{\pag\unhbox\mathbox
\footnote{}{The paper is in final form and no version 
of it will be published elsewhere.} 
\vglue7cc
\box\tbox
\box\abbox
\vskip8pt}

\long\def\abstract#1{{\nine{\bf Abstract.} 
#1

}}

\def\section#1{\vskip-\lastskip\vskip12pt plus2pt minus2pt
{\bf #1}}

\long\def\th#1#2#3{\vskip-\lastskip\vskip4pt plus2pt
{\sc #1} #2\hskip-\lastskip\ {\it #3}\vskip-\lastskip\vskip4pt plus2pt}

\long\def\defin#1#2{\vskip-\lastskip\vskip4pt plus2pt
{\sc #1} #2 \vskip-\lastskip\vskip4pt plus2pt}

\long\def\remar#1#2{\vskip-\lastskip\vskip4pt plus2pt
\sp{#1} #2\vskip-\lastskip\vskip4pt plus2pt}

\def\Proof{\vskip-\lastskip\vskip4pt plus2pt 
\sp{Proo{f.}\ }\ignorespaces}

\def\endproof{\nobreak\kern5pt\nobreak\vrule height4pt width4pt depth0pt
\vskip4pt plus2pt}

\newbox\refbox
\newdimen\refwidth
\long\def\references#1#2{{\nine
\setbox\refbox=\hbox{\nine[#1]}\refwidth\wd\refbox\advance\refwidth by 12pt%
\def\textindent##1{\indent\llap{##1\hskip12pt}\ignorespaces}
\vskip24pt plus4pt minus4pt
\centerline{\bf References}
\vskip12pt plus2pt minus2pt
\parindent=\refwidth
#2

}}

\def\footnoterule{\kern -3pt \hrule width 4cc \kern 2.6pt}

\catcode`@=11
\def\vfootnote#1%
{\insert\footins\bgroup\nine\interlinepenalty\interfootnotelinepenalty%
\splittopskip\ht\strutbox\splitmaxdepth\dp\strutbox\floatingpenalty\@MM%
\leftskip\z@skip\rightskip\z@skip\spaceskip\z@skip\xspaceskip\z@skip%
\textindent{#1}\footstrut\futurelet\next\fo@t}
\catcode`@=12

\def\sthree{{{\cal S}^3}}
\def\tr{{\rm Tr}}
\def\half{{1\over2}}
\def\hol#1{{\rm Hol}_{K_{#1}}(A)}
\def\cd{{\cal D}}
\def\co{{\cal O}}
\def\wilson#1{W_{n_{#1}}^{K_{#1}}(A)}
\def\prodin{\prod_{i=1}^N}
\def\partition{Z(\sthree,L;g_1,g_2,\dots,g_N)}
\def\avewilson#1#2{\left<W_{#1}^{#2}(A)\right>}
\def\one{{\bf 1}}
\def\Keq#1#2{\Bigl\langle\cdots #1_{K_1}\cdots #2_{K_2}\cdots\Bigr\rangle
=\Bigl\langle\cdots #1_{K_1\#_b \tilde K_2}\cdots #2_{K_2}\cdots\Bigr\rangle}
\def\lk{\ell\!k}
\def\M{{\cal M}}
\def\N{{\cal N}}

\def\plusskein{
\setbox0=\hbox{
\setbox1=\hbox{\bigg/}
\setbox2=\hbox{$\backslash$}
\hbox to -\wd1{} \copy1
\kern-\wd1\raise.5\ht1\copy2
\lower\dp1\hbox{\raise\dp2\copy2}
}\copy0}
\def\minusskein{
\setbox0=\hbox{
\setbox1=\hbox{$\bigg\backslash$}
\setbox2=\hbox{/}
\hbox to -\wd1{} \copy1
\kern-\wd1\lower.5\ht1\copy2
\raise\dp1\hbox{\lower\dp2\copy2}
}\copy0}
\def\zeroskein{
\;\bigg\vert\;\bigg\vert\;}
\def\nullskein#1{
\;\bigg\vert #1\;}
\def\clap#1#2{
\setbox0=\hbox{$\displaystyle#2$}
\hbox to \wd0 {\hss$#1$\hss}
\kern-\wd0\hbox to\wd0{\hss$\displaystyle#2$\hss}}


\mathclass{Primary 57R65; Secondary 81T13.}

\abbrevauthors{B.\ Broda}

\abbrevtitle{3-and 4-manifold invariants}

\title{A gauge-field approach\\
to 3- and 4-manifold invariants}

\author{BOGUS\L AW\ BRODA}

\address{Department of Theoretical Physics, University of \L\'od\'z\\
Pomorska 149/153, PL--90-236 \L\'od\'z, Poland\\
E-mail: bobroda@mvii.uni.lodz.pl}

\maketitlebcp

\footnote{}{Research supported by the KBN grants
2P30213906, 2P30221706p01, and the University of
\L\'od\'z grant 505/457.}

\abstract{An approach to construction of topological
invariants of the Reshetikhin-Turaev-Witten type of
3- and 4-dimensional manifolds in the framework of SU(2)
Chern-Simons gauge theory and its hidden (quantum) gauge
symmetry is presented.}

\section{1.\ Intoduction.}
The issue of topological classification of low-dimensional
manifolds, especially of dimensions 3 and 4 (the most
difficult and interesting ones), is a
challenging problem in modern mathematics. One of the most
spectacular events in topology of 3-dimensional manifolds
took place a few years ago, when a new (numerical)
topological invariant of closed orientable
3-manifolds, parametrized by the integer $k$,
defined via {\it surgery} on a framed link, was discovered.
The idea is due to a physicist, Edward Witten, who proposed
the invariant in his famous paper on quantum field theory
and the Jones polynomial [Wit1].
The first explicit and rigorous construction is due to
mathematicians, Reshetikhin and Turaev [RT].
Their approach is combinatorial, whereas non-combinatorial
possibilities, very straightforward though mathematically
less rigorous, are offered by {\it topological} quantum
field theory. The 3-dimensional invariant, known as the {\it
Reshetikhin-Turaev-Witten} (RTW) invariant, is also
frequently referred to as the SU(2){\it -invariant}
because the {\it Kauffman bracket\/} it bases upon (denoted
in mathematical literature as `$\left<\;\right>$') formally
corresponds, in Witten's approach, to the average with
respect to the connection $A$ (defined on a trivial
SU(2) bundle on the 3-manifold $M$) modulo gauge
transformations, weighted by $\exp\left[ikS_{\rm
CS}(A)\right]$. Here $S_{\rm CS}(A)$ is the Chern-Simons
{\it secondary characteristic class}.  Incidentally, the
average is also denoted as `$\left<\;\right>$'. The
construction of the RTW invariant makes use of the {\it
fundamental theorem of surgery} of Lickorish and Wallace on
presentation of every closed connected orientable
3-manifold $M$ via surgery on a framed link in $\sthree$,
and the {\it linear skein theory} associated with the
Kauffman bracket. The construction of the RTW invariant
amounts to showing its invariance with respect to the
{\it Kirby moves} (${\rm K}_1$ and  ${\rm K}_2$).
The Kirby moves are the allowable moves on a framed link,
changing in general (isotopy class of) the link but not
changing the 3-manifold $M$ obtained from the link.
In fact, ${\rm K}_2$ does not even change the 4-manifold
$\cal M$ ($M=\partial\cal M$), because it corresponds to sliding
handles, whereas ${\rm K}_1$ adds a complex projective
space, which changes $\cal M$ but still it does not
change the 3-dimensional boundary $M$. The Kirby move
${\rm K}_1$ means that we can add an unknotted, unlinked
component with framing $\pm 1$, i.e. $\bigcirc_{\pm1},$
whereas ${\rm K}_2$ ammounting to sliding a (upper) line
over a (lower) trivial unknot in an annulus (a bordered
exterior of ``$\bullet$''), which is, in general,
non-trivially immersed in ${\cal S}^3$ is schematically
depicted as
$$
\overline{\odot}\longleftrightarrow
\clap{\odot}{\bigcup}.
$$
\bigskip
In Section 2 we aim to propose a new, heuristic,
non-combinatorial derivation of the RTW invariant in
the framework of non-perturbative (topological) quantum
Chern-Simons (CS) gauge theory [Bro3]. The idea is extremely
simple, and in principle it applies to an arbitrary compact
(semi-)simple Lie group $G$ (not only to the SU(2) one).
Our invariant is essentially the partition function of CS
theory on the 3-manifold $M_L$, defined via {\it surgery} on
the framed link $L$ in the 3-dimensional sphere
$\sthree$.  Actually, surgery instructions are implemented
in the most direct and literal way. The method of cutting
and pasting back is explicitly used
in the standard field-theoretical fashion. Roughly
speaking, cutting corresponds to fixing, whereas pasting
back to identification and summing up the boundary
conditions.

Section 3 is devoted to a 4-dimensional generalization of
the RTW invariant, and is mathematically more rigorous
[Bro1]. In dimension 4, there is a celebrated theorem of
Freedman on classification of closed orientable {\it
simply-connected\/} 4-manifolds, provided by the {\it
intersection form} $Q(\M)$ (and the Kirby-Siebenman
invariant $\alpha(\M)$). The intersection form $Q(\M)$
corresponds to, and for 4-manifolds with boundary defined
via surgery on a link in $\sthree$ is equal to, the {\it
linking matrix} $\lk$. The elements of the symmetric matrix
$\lk$, the {\it linking numbers} (with framing numbers on the
diagonal), are the simplest numerical invariants of a link.
Therefore, one can ask the following questions. Can one use
`non-abelian' invariants of links, for example the Kauffman
bracket polynomial, to obtain `non-abelian' invariants of
4-dimensional manifolds?  Can one extend the idea of RTW to
the 4-dimensional case?  Can one treat simply-connected and
non-simply-connected manifolds uniquely? The answer to
these questions seems to be affirmative. Namely, we propose
an invariant of closed connected orientable
4-manifolds, defined via surgery on a special link in
$\sthree$.  Thus, we have found a quantity invariant with
respect to the 4-dimensional version of the `Kirby moves'.
The idea as well as the construction resembles the original
one, proposed by RTW in the 3-dimensional case, whereas the
4-dimensional version of the {\it Kirby calculus} we need
has been developed by C\'esar~de~S\'a [CdS].

In Section 4, we mention some other known invariants of 3-
and 4-dimensional manifolds.

\section{2.\ Gauge-field approach to 3-manifold
invariants.} 
Our principal goal is to compute the partition function
$Z(M_L)$ of CS theory on the manifold $M_L$, defined via
honest/integer surgery on the framed link
$L=\bigcup_{i=1}^NK_i$ in $\sthree$ (i.e. attaching
2-handles to a 4-ball along a framed link---a particular
case of rational surgery manipulating only tori),
for the SU(2) (gauge) Lie group. Obviously, the
starting point is the partition function [Wit1] of CS theory
$Z(\sthree)$ on the sphere $\sthree$
$$
Z(\sthree)=\int e^{ikS_{\rm CS}(A)}{\cal D}A,
$$
where the functional integration is performed with respect
to the connections $A$ modulo gauge transformations,
defined on a trivial SU(2) bundle on $\sthree$. The classical action is
the {\it CS secondary characteristic class}
$$
S_{\rm CS}(A)={1\over4\pi}\int_\sthree\tr\left(AdA
+{2\over3}A^3\right),
$$
and the expectation value of an observable $\co$ is defined
as 
$$
\left<\co\right>=\int\co
e^{ikS_{\rm CS}(A)}\cd A.
$$

According to the surgery prescription [Rol], we should cut out a
closed tubular neighborhood $N_i$ of $K_i$ (a solid
torus), and paste back a copy of a solid torus $T$,
matching the meridian of $T$ to the (twisted by framing
number) longitude on the boundary torus $\partial N_i$ in
$\sthree$. To this end, in the first step, we should
fix boundary conditions for the field $A$ on the twisted
longitude represented by $K_i$.  Since the only
gauge-invariant (modulo conjugation) quantity defined on a
closed curve is holonomy, we associate the {\it
holonomy} operator $\hol{i}$ to each knot $K_i$.  Thus the
symbol
$$
\partition
$$
should be understood as the {\it constrained} partition
function of CS theory, i.e.\ the values of holonomies along
$K_i$ are fixed
$$
\hol{i}=g_i,\qquad i=1,2,\dots,N.
$$
Now, we can put
$$
\partition
=\left<\prodin\delta(g_i,\hol{i})\right>,
\leqno(1)
$$
where $\delta$ is a (group-theoretic) Dirac
delta-function. Its explicit form following from the
(group-theoretic) Fourier expansion is
$$
\delta(g,h)=\sum_n\overline{\chi_n(g)}\chi_n(h),
\leqno(2)
$$
where $n$ labels inequivalent irreducible
representations (irrep's) of SU(2), and $\chi_n$ is a
character.
Physical observables being used in CS theory are typically
Wilson loops, defined as
$$
\wilson{}=\tr_n(\hol{})\equiv\chi_n(\hol{}).
\leqno(3)
$$
By virtue of (2--3)
$$
\delta(g_i,\hol{i})=\sum_n\overline{\chi_n(g_i)}W_n^{K_i}(A).
\leqno(4)
$$
Inserting (4) into (1) yields, as a basic building
block, the following representation of the constrained
partition function
$$
\partition=\left<\prodin\sum_{n_i}
\overline{\chi_{n_i}(g_i)}\wilson{i}\right>.
$$

In the second step of our construction, we should paste
back the tori matching the pairs of ``longitudes'' (the
twisted longitudes and the meridians), i.~e. we should
identify and sum up the boundary conditions. Since the
interior of a solid torus is homeomorphic to $\sthree$ with
a removed solid torus, actually the meridians play the role
of longitudes in analogous cutting procedures for an unknot
$\left\{\bigcirc\right\}$ (with reversed orientation).
Thus the partition function of CS theory on $M_L$ is
$$
Z(M_L)=\int\prodin dg_i Z(\sthree,\bigcirc;g_i^{-1})
\partition\qquad\qquad\qquad
$$
$$
\qquad\qquad\qquad
=\int\prodin dg_i\sum_{m_i}\sum_{n_i}
\overline{\chi_{m_i}(g_i^{-1})} \,
\overline{\chi_{n_i}(g_i)}
\avewilson{m_i}{\bigcirc}
\left<\prod_{j=1}^N W_{n_j}^{K_j}(A)\right>,
$$
where the reversed orientation of the unknots
$\left\{\bigcirc\right\}$ (corresponding to the meridians
of the pasted back tori) accounts for the power $-1$ of the
group elements $g_i$. From the orthogonality relations for
characters and unitarity of irrep's, it follows that the
3-manifold invariant is of the form
$$
Z(M_L)=\left<\prodin\omega_{K_i}(A)\right>,
\leqno(5)
$$
where
$$
\omega_{K_i}(A)\equiv\sum_{n_i}\avewilson{n_i}{\bigcirc}
\wilson{i}
$$
is an element of the linear skein of an annulus, immersed
in the plane as a regular neighborhood of $K_i$. $\langle
W_n^\bigcirc(A)\rangle$ are some computable coefficients
depending on $n$ and $k$. Eq.~(5) can be
easily generalized to accommodate an ordinary link ${\cal
L}=\bigcup_{i=1}^M{\cal K}_i$ embedded in $M_L$
$$
\left<\prod_{i=1}^M W_{n_i}^{{\cal K}_i}(A)\right>_{M_L}
=\left<\prod_{i=1}^M W_{n_i}^{{\cal K}_i}(A)
\prod_{j=1}^N\omega_{K_j}(A)\right>.
$$

It appears that a very convenient way of organization of
irrep's of SU(2) group is provided by the polynomials
$S_n(x)$, closely related to the Chebyshev polynomials.
$S_n(x)$ are defined recursively by the formula
$$
S_{n+2}(x)=xS_{n+1}-S_n(x),\qquad n=0,1,\dots,
\leqno(6a)
$$
together with the initial conditions
$$
S_0(x)=1,\qquad S_1(x)=x.
\leqno(6b)
$$
By virtue of (6), $S_n(x)$ expresses $n$-th irrep of
SU(2) in terms of powers of the fundamental representation
$x$, denoted as $\one$ henceforth. The explicit solution of (6) is
$$
S_n(2\cos\alpha)={\sin((n+1)\alpha)\over\sin\alpha}.
$$
The skein relations for the fundamental representation
($n=1$) of SU(2) are given by the expression
$$
q^{1\over4}\left<\cdots\plusskein\cdots\right>
-q^{-{1\over4}}\left<\cdots\minusskein\cdots\right>
=(q^{\half}-q^{-\half})
\left<\cdots\zeroskein\cdots\right>,
\leqno(7a)
$$
$$
\left<\cdots\nullskein{\pm1}\cdots\right>
=-q^{\pm{3\over4}}\left<\cdots\nullskein{0}\cdots\right>,
\leqno(7b)
$$
where the integers in (7b) mean framings,
and $q=\exp{2\pi i\over k}$. Closing 
the left legs of all the (three) diagrams in (7a) with
arcs, as well as the right ones, next applying (7b), and
using the property of locality, we obtain
$$
-(q-q^{-1})\avewilson{\one}{\bigcirc}
=(q^{\half}-q^{-\half})\avewilson{\one}{\bigcirc\bigcirc}
=(q^{\half}-q^{-\half})\avewilson{\one}{\bigcirc}^2.
$$
Hence
$$
\avewilson{\one}{\bigcirc}=-\left(q^{\half}+q^{-\half}\right)
=-2\cos{\pi\over k},
$$
and by virtue of the so-called satellite formula
$$
\avewilson{n}{\bigcirc}=S_n\left(-2\cos{\pi\over k}\right)
=(-)^n{\sin{(n+1)\pi\over k}\over\sin{\pi\over k}}
=(-)^n{q^{n+1\over2}-q^{-{n+1\over2}}\over
q^{\half}-q^{-\half}}.
\leqno(8)
$$
We can observe a remarkable property of (8) for
$n=k-1$, namely
$$
\avewilson{k-1}{\bigcirc}=0.
\leqno(9)
$$
It appears that for any $\cal K$
$$
\langle\cdots W_{k-1}^{\cal K}(A)\cdots\rangle=0.
\leqno(10)
$$
In particular, Eq.~(10) immediately follows from (9)
for any $\cal K$ that can be unknotted with corresponding
skein relations.  Thus we can truncate representations
of SU(2) above the value $k-2$, and assume
$$
0\leq n\leq k-2,\qquad k=2,3,\dots.
$$
The final explicit form of $\omega_K$ for
the group SU(2) is then
$$
\omega_K(A)=\sum_{n=0}^{k-2}(-)^n
{q^{n+1\over2}-q^{-{n+1\over2}}\over
q^{\half}-q^{-\half}}S_n\left(W_\one^K(A)\right).
\leqno(11)
$$
Strictly speaking,
$Z(M_L)$ is invariant with respect to the second Kirby
move ${\rm K}_2$. It means that it is insensitive to the operation of
sliding one of its handles over another one. But up to now we
have not considered the issue of the determination of
normalization. It appears that proper normalization of the
partition function $Z(M_L)$ universally follows from the
requirement of its invariance with respect to the first
Kirby move ${\rm K}_1$.

The approach proposed above differs from Witten's one
[Wit1] in that we explicitly construct the invariant
via the cutting and pasting procedure using standard
field-theoretic tools.

\section{3.\ Generalization to 4 dimensions.}
An arbitrary closed connected orientable
4-manifold $\M$ can be obtained via surgery in $\sthree$
on a {\it special framed link} $\left(L,f\right)$ [CdS].

\defin{Definition.}{The special framed link $L$ is a sum of
two sorts of knots 
$$
L=\bigcup_{i=1}^N K_i
\cup
\bigsqcup_{i=1}^{\dot N} \dot K_i,
$$
where $\bigl\{K_i\bigr\}_{i=1}^N$ are ordinary knots (corresponding to 2-handles), and
$\bigl\{\dot K_i\bigr\}_{i=1}^{\dot N}$ are special knots (corresponding to 1-handles).
The special knots, denoted with dotes, are trivial (with
zero framing), and mutually unlinked unknots, and the whole
link, when regarded as a description of a 3-manifold,
represents a connected sum of copies of
${\cal S}^1\times {\cal S}^2$. The symbol ``$\sqcup$"
means the ``distant sum"---the components are mutually
unlinked.}

Now, we introduce the following decomposition ({\it
gradation}) of $\omega$, defined in (11), into an even
($+$) and odd ($-$) parts (integer and half-integer `spins', respectively)
$$
\omega=\omega^+ + \omega^-.
$$

Let us denote as $\tilde K$ the result of pushing a knot $K$
off itself (missing the rest of the link $L$) using the
framing $f$ of $K$, whereas as $K_1\#_b K_2$ a (band)
connected sum of the two knots $K_1$, $K_2$, where $b$ is
any band missing the rest of $L$.

\th{Proposition.}{}{For arbitrary complex numbers, $a^+$, $a^-$, we have
the following `Kirby calculus'
$$
\Keq{\alpha}{\omega},
$$
$$
\Keq{\alpha^2}{\left(a^+ \omega^+ + a^-\omega^-\right)},
$$
where $\alpha$ is a linear combination of $S_n$ with
arbitrary coefficients (compare with the definition of $\omega$ (11)).}

\remar{Remark.}{The first equality expresses a standard property of
$\omega$, whereas the second one follows from the
observation that the even element $\alpha^2$ respects the
gradation in all cablings.}

\th{Corollary.}{}{From the Proposition we can derive the following
`Kirby equalities'
$$
\Keq{\omega}{\omega},
$$
$$
\Keq{\omega^+}{\omega^+},
$$
$$
\Keq{\omega^+}{\omega}.
$$
}
The fourth equality,
$$
\Keq{\omega}{\omega^+},
$$
is, in general, not true.

Henceforth, $H$ and $\dot H$ are two components of the {\it
special Hopf\/} link $\cal H$, ordinary and special one
respectively, ${\cal H}=H\cup\dot H$.  

\th{Theorem.}{Let $\nu$ be the nullity of the (extended)
linking matrix $\lk$. Then
$$
I_k(\M)=
{\left<\prod_{i=1}^N \omega_{K_i}^+
\prod_{i=1}^{\dot N} \omega_{\dot K_i}\right>
\over
\left< \omega_\bigcirc^+ \right>^\nu
\left< \omega_H^+ \omega_{\dot H} \right>^{(N+\dot N-\nu)/2}}
$$
is an invariant of the closed, connected, orientable
4-manifold $\M=\M_L$, a complex number parametrized by
the integer $k$, independent of the choice of the
representative $\left(L,f\right)$.}

Below, we give a list of all the allowable
`4-dimensional Kirby moves', so-called $\Gamma$-{\it
moves} [CdS]: 

($a$)
sliding one of the special knots over another special one; 

($b$)
sliding one of the ordinary knots over one of the special ones;

($c$)
sliding one of the ordinary knots over another ordinary
one; 

($d$)
introducing or deleting a special Hopf link;

($e$)
introducing or deleting a trivial unknot;

($f$)
isotoping the link picture in $\sthree$.

\Proof
We should show that $I_k(\M)$ is invariant with respect to all
the $\Gamma$-moves. $a$-, $b$- and $c$-invariance of $I_k(\M)$
immediately follows from the Corollary as well as from the
invariance of $N$, $\dot N$ and $\nu$. $d$-invariance is a
consequence of the following transformation rule of the
linking matrix $\lk$, accompanying the introduction of a
special Hopf link $\cal H$,
$$
\lk\longrightarrow\left(\matrix{\lk&0&0\cr
                                  0&0&1\cr
                                  0&1&0\cr}\right).
$$
Hence the corresponding shift of the dimension and nullity
of $\lk$
$$
N\longrightarrow N+1
$$
$$
\dot N\longrightarrow\dot N+1
$$
$$
\nu\longrightarrow\nu,
$$
compensates the (factorized out) Kauffman bracket in the
numerator. Similarly, $e$-invar\-iance corresponds to the
transformation rule
$$
\lk\longrightarrow\left(\matrix{\lk&0\cr
                                  0&0\cr}\right),
$$
and consequently the shift
$$
N\longrightarrow N+1
$$
$$
\dot N\longrightarrow\dot N
$$
$$
\nu\longrightarrow\nu+1,
$$
also compensates the numerator. $f$-invariance directly
follows from fundamental properties of the Kauffman bracket and
the linking matrix $\lk$.
\vskip4pt plus2pt

\remar{Remark 1.}{The invariant $I_k(\M)$ possesses the following obvious
properties: 

\noindent
(1)\quad
{\it Multiplicativity},
$$
I_k(\M\#\N)=I_k(\M)\cdot I_k(\N),
$$
(2)\quad
{\it Orientation sensitivity},
$$
I_k(\overline{\M})=\overline{I_k(\M)},
$$
(3)\quad
{\it Normalization},
$$
I_k({\cal S}^4)=1,
$$
where $\M\#\N$ denotes a connected sum of $\M$ and $\N$.}

\remar{Remark 2.}{Crane and Yetter [CY] have found a
4-dimensional topological invariant, defined via
triangulation, which is basically equivalent to ours.}

\remar{Remark 3.}{It has been proved that $I_k(\M)$ is
expressible by classical invariants (signature and Euler
character) [CKY2], but there is a gap in the proof announced in [CKY1].}

\section{4.\ Final remarks.}
In the 3-dimesnional case, besides the original RTW
invariant, we have some other invariants of the same CS
origin, and therefore more or less mutually related.
The so-called {\it Turaev-Viro} invariant [TV] (the
Crane-Yetter invariant, mentioned in Remark 2,
is its 4-dimensional counterpart) can be calculated
from triangulation of the manifold, and corresponds
to the square of the modulus of the RTW invariant.
Whereas the Kohno invariant can be calculated from
the Heegaard decomposition of the manifold, and is
also basically equivalent to the RTW invariant [Bro2].
Perturbative expansion provides us with some further
family of invariants, defined for homology spheres,
so-called {\it Ohtsuki-Garoufalidis} invariants [Oht].
The first perturbative term of this family is the famous
{\it Casson(-Walker)} invariant, which was originally
defined non-perturbatively, and has a surgical description.

In the 4-dimensional situation, we have also a distinct
world of `differentiable invariants', the famous {\it
Donaldson} invariants and their  building block,
the {\it Seiberg-Witten} invariant [Wit2],
mathematically described in this volume.
But up to now, it is no clear whether
the `differentiable invariants' are related to the combinatorial ideas presented here.


\references{CKY2}
{

\item{[Bro1]} B. \spa{Broda},
{\it A surgical invariant of 4-manifolds,
in: Proc. Conf. on Quantum Topology, Kansas 1993}\/,
D.~N.~Yetter (ed.), World Scientific, Singapore, 1994, 45--50. 

\item{[Bro2]} B. \spa{Broda},
{\it TQFT versus RCFT: 3-D topological invariants}\/,
Mod. Phys. Lett. A 10 (1995), 331--336.

\item{[Bro3]} B. \spa{Broda},
{\it Chern-Simons approach to three-manifold invariants}\/, 
Mod. Phys. Lett. A 10 (1995), 487--493.

\item{[CdS]} E. \spa{C{\'e}sar~ de~ S{\'a}},
{\it A link calculus for 4-manifolds,
in: Topology of low-dimensional manifolds, Proc. Sec. Conf.
Sussex}\/, Lecture Notes in Math. 722, Springer, Berlin, 1979, 16--30. 

\item{[CKY1]} L. \spa{Crane}, L. H. \spa{Kauffman} and D. \spa{Yetter},
{\it $U_q(sl_2)$ Invariants at Principal and Non-principal Roots of Unity}\/,
Adv. in App. Clifford Algebras 3 (1993), 223.

\item{[CKY2]} L. \spa{Crane}, L. H. \spa{Kauffman} and D. \spa{Yetter},
{\it On the Classicality of Broda's SU(2) Invariants of 4-Manifolds}\/,
Adv. in App. Clifford Algebras 3 (1993), 223.

\item{[CY]} L. \spa{Crane} and D. \spa{Yetter},
{\it A categorical construction of 4d topological quantum
field theories, in: Quantum Topology}\/,
L. Kauffman and R. Baadhio (eds.), World Scientific, Singapore, 1993.

\item{[Oht]} T. \spa{Ohtsuki},
{\it A polynomial invariant of integral homology 3-spheres}\/,
Math. Proc. Camb. Phil. Soc. 117 (1995), 83--112.

\item{[Rol]} D. \spa{Rolfsen},
{\it Knots and Links}\/,
Publish or Perish, Wilmington, 1976, Chapt.~9.

\item{[RT]} N. \spa{Reshetikhin} and V. G. \spa{Turaev},
{\it Invariants of 3-manifolds via link polynomials and quantum groups}\/,
Invent. math. 103 (1991), 547--597.

\item{[Wit1]} E. \spa{Witten},
{\it Quantum Field Theory and the Jones polynomial}\/,
Commun. Math. Phys. 121 (1989), 351--399.

\item{[Wit2]} E. \spa{Witten},
{\it Monopoles and four-manifolds}\/,
Math. Res. Lett. 1 (1994), 769.

\item{[TV]} V. G. \spa{Turaev} and O. Y. \spa{Viro},
{\it State sum invariants of 3-manifolds and quantum 6j-symbols}\/,
Topology 31 (1992), 865--902.

}

\bye